\documentclass[a4, 10pt]{amsart}
\usepackage{amssymb}
\usepackage{amstext}
\usepackage{amsmath}
\usepackage{amscd}
\usepackage{latexsym}
\usepackage{amsfonts}

\theoremstyle{plain}
\newtheorem{thm}{Theorem}[section]
\newtheorem*{thm*}{Theorem}
\newtheorem*{cor*}{Corollary}

\newtheorem*{thma}{Theorem A}

\newtheorem*{thmb}{Theorem B}

\newtheorem{prop}[thm]{Proposition}
\newtheorem{lem}[thm]{Lemma}

\newtheorem{cor}[thm]{Corollary}
\newtheorem{claim}{Claim}
\newtheorem*{claim*}{Claim}

\theoremstyle{definition}
\newtheorem{defn}[thm]{Definition}

\newtheorem{rem}[thm]{Remark}

\theoremstyle{remark}
\newtheorem*{pf}{{\sl Proof}}

\newtheorem*{cpf}{{\sl Proof of Claim}}

\numberwithin{equation}{thm}
\def\Hom{\mathrm{Hom}}

\def\Ext{\mathrm{Ext}}

\def\End{\mathrm{End}}

\def\mod{\mathrm{mod}}
\def\Coker{\mathrm{Coker}}

\def\Im{\mathrm{Im}}
\def\tr{\mathrm{Tr}}

\def\m{\mathfrak m}

\def\p{\mathfrak p}

\def\Gdim{\mathrm{G}\mathrm{dim}}
\def\GCdim{\mathrm{G}_C\mathrm{dim}}
\def\GCpdim{\mathrm{G}_{C_{\p}}\mathrm{dim}}
\def\Cdim{C\mathrm{dim}}
\def\Cpdim{C_\p\mathrm{dim}}
\def\odim{\omega\mathrm{dim}}
\def\Rdim{R\mathrm{dim}}

\def\Kdim{\mathrm{dim}}
\def\depth{\mathrm{depth}}
\def\Supp{\mathrm{Supp}}

\def\Ass{\mathrm{Ass}}

\def\pd{\mathrm{pd}}
\def\id{\mathrm{id}}

\def\height{\mathrm{ht}}
\def\grade{\mathrm{grade}}

\def\Spec{\mathrm{Spec}}
\def\add{\mathrm{add}}

\def\X{{\mathcal X}}
\def\Y{{\mathcal Y}}

\def\xx{\text{\boldmath $x$}}

\begin{document}


\title[Spherical and Cohen-Macaulay approximations]{A new approximation theory which unifies\\
spherical and Cohen-Macaulay approximations}
\author{Ryo Takahashi}
\address{Department of Mathematics, School of Science and Technology, Meiji University, Kawasaki 214-8571, Japan}
\email{takahasi@math.meiji.ac.jp}
\thanks{
\endgraf
}
\maketitle
\begin{abstract}
This paper gives a new approximation theory for finitely generated modules over commutative Noetherian rings, which unifies two famous approximation theorems; one is due to Auslander and Bridger and the other is due to Auslander and Buchweitz. Modules admitting such approximations shall be studied.
\end{abstract}
\section{Introduction}

In the late 1960s Auslander and Bridger \cite{ABr} introduced a notion of approximation which they used to prove that every module whose $n^{\rm th}$ syzygy is $n$-torsionfree can be described as the quotient of an $n$-spherical module by a submodule of projective dimension less than $n$. About two decades later Auslander and Buchweitz \cite{ABu} introduced the notion of Cohen-Macaulay approximation which they used to show that the category of finitely generated modules over a Cohen-Macaulay local ring with the canonical module is obtained by gluing together the subcategory of maximal Cohen-Macaulay modules and the subcategory of modules of finite injective dimension. Our purpose is to give a new approximation theorem, which unifies these two notions of \cite{ABr} and \cite{ABu}. Before stating our own result, let us briefly summarize the theorems of \cite{ABr} and \cite{ABu}.

Throughout let $R$ be a commutative Noetherian ring. Let $M$ be a finitely generated $R$-module and $n \geq 1$ an integer. Then we say that $M$ is $n$-spherical if $\Ext _R ^i(M,R)=0$ for $1\leq i\leq n$. We say that $M$ is $n$-torsionfree if the transpose $\tr M$ of $M$ is $n$-spherical. Let $\Omega ^nM$ denote the $n^{\rm th}$ syzygy of $M$. With this notation the approximation theorem of Auslander and Bridger is stated as follows.

\begin{thm}\cite{ABr}\label{abr}
The following are equivalent for a finitely generated $R$-module $M$.
\begin{enumerate}
\item[{\rm (1)}]
$\Omega ^nM$ is $n$-torsionfree.
\item[{\rm (2)}]
There exists an exact sequence $0 \to Y \to X \to M \to 0$ of $R$-modules such that $X$ is $n$-spherical and $\pd _R\,Y<n$.
\end{enumerate}
\end{thm}

\noindent
Let us call an exact sequence as above an $n$-spherical approximation of $M$.

The Cohen-Macaulay Approximation Theorem of Auslander and Buchweitz is the following.

\begin{thm}\cite{ABu}\label{abu}
Let $R$ be a Cohen-Macaulay local ring with the canonical module.
Then for every finitely generated $R$-module $M$ there exists an exact sequence $0 \to Y \to X \to M \to 0$ of $R$-modules such that $X$ is a maximal Cohen-Macaulay $R$-module and $\id _R\,Y<\infty$.
\end{thm}

\noindent
We call an exact sequence as above a Cohen-Macaulay approximation of $M$.

For a moment suppose that $R$ is a Gorenstein local ring of dimension $d$. Then a finitely generated $R$-module has finite projective dimension if and only if it has finite injective dimension. The $d^{\rm th}$ syzygy of any finitely generated $R$-module is $d$-torsionfree because it is a maximal Cohen-Macaulay $R$-module.
The local duality theorem guarantees that a finitely generated $R$-module is a maximal Cohen-Macaulay $R$-module if and only if it is $d$-spherical.
Thus, although in general Theorems \ref{abr} and \ref{abu} have no implications to each other, Theorem \ref{abr} implies Theorem \ref{abu} in the case where $R$ is a Gorenstein local ring, so that they yield the common consequence which asserts the existence of Cohen-Macaulay approximations over Gorenstein local rings. Added to it, there is a strong similarity between the two exact sequences in Theorems \ref{abr} and \ref{abu}. It seems natural to guess that behind the above two theorems there is hidden a common source unifying them, which we are going to report in this paper.

Let us now state our own results, explaining how this paper is organized. In Section 2 we will make some definitions and state basic properties, which we need throughout this paper.

In Section 3 we will prove our main result. To state the result precisely we need new notation. Let $M$ and $C$ be finitely generated $R$-modules and let $n \geq 1$ be an integer. Let $\lambda _M$ be the natural homomorphism $M\to \Hom _R (\Hom _R(M,C),C)$. Then we say that $M$ is $1$-$C$-torsionfree if $\lambda _M$ is a monomorphism, and that $M$ is $n$-$C$-torsionfree if $n\geq 2$, $\lambda _M$ is an isomorphism, and $\Ext _R^i(\Hom _R(M,C), C)=0$ for $1\leq i\leq n-2$. The $R$-module $C$ is called $1$-semidualizing if $\lambda _R$ is a monomorphism and $\Ext _R ^1(C,C)=0$, and $C$ is called $n$-semidualizing if $n\geq 2$, $\lambda _R$ is an isomorphism, and $\Ext _R ^i(C,C)=0$ for $1\leq i\leq n$. The $R$-module $M$ is said to be $n$-$C$-spherical if $\Ext _R ^i(M,C)=0$ for $1\leq i\leq n$. We denote by $\Cdim _R\,M$ the infimum of nonnegative integers $m$ such that there exists an exact sequence $0 \to C_m \to C_{m-1} \to \cdots \to C_0 \to M \to 0$ of $R$-modules, where each $C_i$ is a direct summand of a finite direct sum of copies of $C$. The main result of this paper, which we shall prove in Section 3 (Theorem \ref{main}), is stated as follows.

\begin{thma}
Let $M$ and $C$ be finitely generated $R$-modules and $n \geq 1$ an integer. Assume that $C$ is $n$-semidualizing. Then the following two conditions on $M$ are equivalent to each other.
\begin{enumerate}
\item[{\rm (1)}]
$\Omega ^nM$ is $n$-$C$-torsionfree.
\item[{\rm (2)}]
There exists an exact sequence $0 \to Y \to X \overset{f}{\to} M \to 0$ of $R$-modules such that $X$ is $n$-$C$-spherical and $\Cdim _R\,Y<n$.
\end{enumerate}
\end{thma}

\noindent
The above map $f$ is a right approximation over the full subcategory of $\mod\,R$ consisting of all $n$-$C$-spherical $R$-modules, where $\mod\,R$ denotes the category of finitely generated $R$-modules. We call an exact sequence as above an $n$-$C$-spherical approximation of $M$.

Here let us briefly explain how Theorem A implies both Theorems \ref{abr} and \ref{abu}. We notice that $R$ is an $n$-semidualizing $R$-module and that $\Rdim _R\,M=\pd _R\,M$ for a finitely generated $R$-module $M$. Our Theorem A thus implies Theorem \ref{abr}. Let $R$ be a $d$-dimensional Cohen-Macaulay local ring with the canonical module $\omega$ and $M$ a finitely generated $R$-module. Then $\omega$ is $d$-semidualizing, and $\Omega ^dM$ is $d$-$\omega$-torsionfree since it is a maximal Cohen-Macaulay $R$-module. The local duality theorem implies that an $R$-module $X$ is $d$-$\omega$-spherical if and only if $X$ is a maximal Cohen-Macaulay $R$-module. For a finitely generated $R$-module $Y$, $\odim _R\,Y<\infty$ if and only if $\id _R\,Y<\infty$. Our Theorem A thus implies Theorem \ref{abu}.

Sections 4 and 5 of this paper are devoted to the analysis of the problem of when the first condition of Theorem A is satisfied. In Section 4 we shall closely study the structure of finitely generated $R$-modules whose $n^{\rm th}$ syzygies are $n$-$C$-torsionfree. We will give some results on $n$-$C$-torsionfree modules corresponding to those on $n$-torsionfree modules given by \cite{ABr} and \cite{Masek}.

In Section 5 we will prove the following result, which is  the second main result of this paper (Theorem \ref{cond}).

\begin{thmb}
Let $C$ be a finitely generated $R$-module such that $R$ is $n$-$C$-torsionfree. Then the following two conditions on $C$ are equivalent to each other.
\begin{enumerate}
\item[{\rm (1)}]
$\Omega ^nM$ is $n$-$C$-torsionfree for any $R$-module $M$.
\item[{\rm (2)}]
$\id _{R_\p}\,C_\p <\infty$ for any $\p\in\Spec\,R$ with $\depth\,R_\p \leq n-2$.
\end{enumerate}
\end{thmb}

\noindent
Combining Theorems A and B, we see that if $C$ is an $n$-semidualizing $R$-module such that $\id _{R_\p}\,C_\p <\infty$ for any $\p\in\Spec\,R$ with $\depth\,R_\p \leq n-2$, then the $n$-$C$-spherical $R$-modules form a contravariantly finite subcategory of $\mod\,R$.

In what follows, let $R$ denote a commutative Noetherian ring. All $R$-modules considered are finitely generated. 

\section{Basic definitions}

Let $M$ be an $R$-module and let
$$
\cdots \to P_n \overset{d_n}{\to} P_{n-1} \to \cdots \to P_0 \to M \to 0
$$
be a projective resolution of $M$.
Let $\Omega ^n M$ be the image of $d_n$, which we call the {\it $n^{\rm th}$ syzygy} of $M$.
This is uniquely determined up to projective summand.
An $R$-module $N$ is said to be {\it $n$-syzygy} if $N$ is isomorphic to the $n^{\rm th}$ syzygy of some $R$-module $M$.
Let $(-)^{\ast}=\Hom _R (-, R)$ denote the $R$-dual.
We put $\tr M = \Coker (d_1^{\ast} : P_0^{\ast} \to P_1^{\ast})$ and call it the {\it transpose} of $M$.
This is also uniquely determined up to projective summand.
See \cite{ABr} and \cite{Masek} for more details on transposes.

\begin{defn}\cite{ABr}
Let $n \geq 1$ be an integer and $M$ an $R$-module.
Then $M$ is said to be $n$-torsionfree if $\Ext _R ^i(\tr M,R)=0$ for all $1\leq i\leq n$.
\end{defn}

Let $M$ be an $R$-module.
Then there exists an exact sequence $$0 \to \Ext _R ^1(\tr M, R) \to M \overset{\rho _M}{\to} M^{\ast\ast} \to \Ext _R ^2(\tr M, R) \to 0$$ of $R$-modules with $\rho _M$ the natural homomorphism (\cite[Proposition (2.6)]{ABr}).
Therefore, since $M^{\ast}$ is isomorphic to $\Omega ^2(\tr M)$ up to projective summand, $M$ is $1$-torsionfree if and only if $\rho _M$ is a monomorphism.
When $n\geq 2$, $M$ is $n$-torsionfree if and only if $\rho _M$ is an isomorphism and $\Ext _R ^i (M^{\ast}, R)=0$ for all $1\leq i\leq n-2$.

In what follows, unless otherwise specified, let $n \geq 1$ be an integer, $C$ an $R$-module, $(-)^{\dag}=\Hom _R (-, C)$ the $C$-dual, and $\lambda _M : M \to M^{\dag\dag}$ the natural homomorphism ($M$ an $R$-module).
We then generalize the notion of an $n$-torsionfree module.

\begin{defn}
Let $M$ be an $R$-module.\\
(1) We say that $M$ is $1$-$C$-torsionfree if $\lambda _M$ is a monomorphism.\\
(2) Suppose that $n\geq 2$. Then we say that $M$ is $n$-$C$-torsionfree if $\lambda _M$ is an isomorphism and $\Ext _R^i(M^\dag, C)=0$ for all $1\leq i\leq n-2$.
\end{defn}

\noindent
Every $n$-$C$-torsionfree $R$-module is $i$-$C$-torsionfree for all $1\leq i\leq n$. The $n$-$R$-torsionfree property is the same as the $n$-torsionfree property. If $R$ is a Cohen-Macaulay local ring with the canonical module $\omega$, then a maximal Cohen-Macaulay $R$-module is $n$-$\omega$-torsionfree for every $n\geq 1$.

\begin{defn}\label{sd}
(1) We say that $C$ is $1$-semidualizing if $\lambda _R: R \to \Hom_R(C,C)$ is a monomorphism and $\Ext _R ^1(C,C)=0$.\\
(2) Suppose that $n \geq 2$. Then we say that $C$ is $n$-semidualizing if $\lambda _R : R \to \Hom_R(C,C)$ is an isomorphism and $\Ext _R ^i(C,C)=0$ for all $1\leq i\leq n$.
\end{defn}

\noindent
An $n$-semidualizing $R$-module is $i$-semidualizing for all $1\leq i\leq n$.
If $C$ is $n$-semidualizing, then $R$ is $n$-$C$-torsionfree.
The $R$-module $R$ is an $n$-semidualzing $R$-module for every $n\geq 1$.
Recall that an $R$-module $C$ is called {\it semidualizing} if $\lambda _R$ is an isomorphism and $\Ext _R ^i(C,C)=0$ for any $i\geq 1$.
A semidualizing $R$-module is $n$-semidualizing for every $n\geq 1$.
In particular, the canonical module of a Cohen-Macaulay local ring is $n$-semidualizing for every $n\geq 1$.

The following proposition says that there are a lot of $n$-semidualizing modules.
Recall that a local ring $(R,\m)$ has an {\it isolated singularity} if $R_\p$ is a regular local ring for any $\p\in\Spec\,R-\{\m\}$.

\begin{prop}\label{goto}
Let $R$ be a Cohen-Macaulay local ring of dimension $d\geq 2$ with an isolated singularity.
Let $I$ be an ideal of $R$ which is a maximal Cohen-Macaulay $R$-module.
Then $\lambda _R : R\to \End _R(I)$ is an isomorphism and $\Ext _R ^i(I,I)=0$ for every $1\leq i\leq d-2$.
Hence $R$ is $d$-$I$-torsionfree, and $I$ is $(d-2)$-semidualizing.
\end{prop}

\begin{pf}
Note that $R$ is a normal domain and that $I$ is a nonzero ideal of $R$.
Hence it follows by \cite[pp. 220--221]{GR} or \cite[Theorem 2.1]{deJong} that $\lambda _R$ is an isomorphism.
Let us prove that $\Ext _R^i(I,I)=0$ for $1\leq i\leq d-2$.
This statement can be shown by using the proof of \cite[Proposition 2.5.1]{Iyama}.
Set $r=\sup\{\,n\,|\,\Ext _R^i(I,I)=0\text{ for }1\leq i\leq n\,\}$.
We want to show that $r\geq d-2$.
Suppose that $r<d-2$.
Take a free resolution $F_{\bullet}$ of the $R$-module $I$.
Dualizing $F_{\bullet}$ by $I$, we obtain an exact sequence
\begin{align*}
0 & \to R \to \Hom _R(F_0,I) \overset{\delta _r}{\to} \cdots \overset{\delta _2}{\to} \Hom _R(F_{r-1},I) \overset{\delta _1}{\to} \Hom _R(F_r,I) \\
& \overset{\delta _0}{\to} \Hom _R(\Omega ^{r+1}I,I) \to \Ext _R^{r+1}(I,I) \to 0
\end{align*}
Put $N_i=\Im\,\delta _i$ for $0\leq i\leq r$.
The definition of $r$ implies that $\Ext _R^{r+1}(I,I)\neq 0$.
Since $R$ has an isolated singularity and $I$ is maximal Cohen-Macaulay, the $R$-module $\Ext _R ^{r+1}(I,I)$ has finite length.
Hence we have $\depth _R\,\Ext _R^{r+1}(I,I)=0$.
As $\depth _R\,\Hom _R(\Omega ^{r+1}I,I)\geq\min\{\, 2,\,\depth _R\,I\,\}=2>0$ (cf. \cite[Exercise 1.4.19]{BH}), we obtain $\depth _R\, N_0=1$ by the depth lemma.
Noting that each $\Hom _R(F_i, I)$ is a maximal Cohen-Macaulay $R$-module, by the depth lemma, we get $\depth _R\,N_i=i+1$ for $0\leq i\leq r$, and $d=\depth\,R=r+2<d$.
This is a contradiction, which shows that $r\geq d-2$ and the proof is completed.
\qed
\end{pf}

We denote by $\mod\,R$ the category of finitely generated $R$-modules.
Let $\X$ be a full subcategory of $\mod\,R$.
An $R$-homomorphism $f:X\to M$ is called a {\it right $\X$-approximation} of $M$ if $X$ belongs to $\X$ and the sequence $\Hom _R (-,X) \overset{(-,f)}{\longrightarrow} \Hom _R (-,M) \to 0$, where $(-,f)=\Hom _R(-,f)$, is exact on $\X$.
We say that $\X$ is {\it contravariantly finite} if any $X\in\X$ has a right $\X$-approximation.
An $R$-complex $(\cdots \overset{f_2}{\to} X_1 \overset{f_1}{\to} X_0 \overset{f_0}{\to} M)$ is called a {\it right $\X$-resolution} of $M$ if each $X_i$ belongs to $\X$ and the sequence $\cdots \overset{(-,f_2)}{\longrightarrow} \Hom _R (-, X_1) \overset{(-,f_1)}{\longrightarrow} \Hom _R (-, X_0) \overset{(-,f_0)}{\longrightarrow} \Hom _R (-, M) \to 0$ is exact on $\X$.
A {\it left $\X$-approximation}, a {\it covariantly finite} subcategory and a {\it left $\X$-resolution} are defined dually.

For an $R$-module $X$, we denote by $\add\,X$ the full subcategory of $\mod\,R$ consisting of all direct summands of finite direct sums of copies of $X$.
For an $R$-module $M$, we define $\Cdim _R\,M$, the {\it $\add\,C$-resolution dimension} of $M$, to be the infimum of nonnegative integers $n$ such that there exists an exact sequence $0 \to C_n \to C_{n-1} \to \cdots \to C_0 \to M \to 0$ with each $C_i$ being in $\add\,C$.
Note that $\add\,R$-resolution dimension is the same as projective dimension.
We make the following definition.

\begin{defn}
Let $M$ be an $R$-module.\\
(1) We say that $M$ is $n$-spherical if $\Ext _R ^i(M,R)=0$ for all $1\leq i\leq n$.
We call an exact sequence $0 \to Y \to X \to M \to 0$ of $R$-modules an $n$-spherical approximation if $X$ is $n$-spherical and $\pd _R\,Y<n$.\\
(2) We say that $M$ is $n$-$C$-spherical if $\Ext _R ^i(M,C)=0$ for all $1\leq i\leq n$.
We call an exact sequence $0 \to Y \to X \to M \to 0$ of $R$-modules an $n$-$C$-spherical approximation if $X$ is $n$-$C$-spherical and $\Cdim _R\,Y<n$.
\end{defn}

\noindent
The $n$-$R$-spherical property is the same as the $n$-spherical property.
By virtue of the local duality theorem, if $R$ is a Cohen-Macaulay local ring with the canonical module $\omega$, then an $R$-module is maximal Cohen-Macaulay if and only if it is $d$-$\omega$-spherical.

\section{The approximation theorem}

In this section, we will discuss when a given module has an $n$-$C$-spherical approximation.
We shall actually give an equivalent condition for a module to have an $n$-$C$-spherical approximation in the case where $C$ is an $n$-semidualizing module.
Using this equivalent condition, we will prove two well-known approximation theorems: one is due to Auslander and Bridger, and the other is due to Auslander and Buchweitz.

First of all, we make a remark, and characterize $n$-$C$-torsionfree modules by using left $\add\,C$-resolutions.

\begin{rem}\label{triple}
The map $\lambda _{M^\dag}: M^{\dag}\to M^{\dag\dag\dag}$ is a split monomorphism for any $R$-module $M$.
Indeed, it is easy to check that the composite map $(\lambda _M)^{\dag}\cdot\lambda _{M^\dag}$ is the identity map of $M^{\dag}$.
\end{rem}

\begin{prop}\label{resol}
Let $M$ be an $R$-module.
Then the following statements hold.
\begin{enumerate}
\item[{\rm (1)}]
$M$ is $1$-$C$-torsionfree if and only if $M$ has an injective left $\add\,C$-approximation.
\item[{\rm (2)}]
\begin{enumerate}
\item[{\rm (i)}]
If $M$ is $2$-$C$-torsionfree, then $M$ has an exact left $\add\,C$-resolution $0 \to M \to C_0 \to C_1$.
\item[{\rm (ii)}]
The converse holds if $\lambda _C$ is an isomorphism.
\end{enumerate}
\item[{\rm (3)}]
Let $n\geq 3$.
Suppose that $\lambda _C$ is an isomorphism and $\Ext _R ^i(C^\dag,C)=0$ for all $1\leq i\leq n-3$.
\begin{enumerate}
\item[{\rm (i)}]
If $M$ is $n$-$C$-torsionfree, then $M$ has an exact left $\add\,C$-resolution $0 \to M \to C_0 \to C_1 \to \cdots \to C_{n-1}$.
\item[{\rm (ii)}]
The converse holds if $\Ext _R ^{n-2}(C^\dag, C)=0$.
\end{enumerate}
\end{enumerate}
\end{prop}

\begin{pf}
(1) Suppose that $M$ is $1$-$C$-torsionfree.
Dualizing a free cover $R^r\to M^{\dag}$ of $M^\dag$ by $C$, we have an injection $\alpha :M^{\dag\dag}\to C^r$.
Hence we get an injection $\beta = \alpha\cdot\lambda _M:M \to C^r$.
Let $f_1, \dots, f_l$ be a system of generators of the $R$-module $\Hom _R(M,C)$.
Taking the direct sum of $f_1, \dots, f_l$, we construct a homomorphism $f:M\to C^l$.
It is easily seen that $f$ is a left $\add\,C$-approximation of $M$.
In particular, the $C$-dual homomorphism $f^{\dag}$ is surjective.
Since $\beta$ factors through $f$, the homomorphism $f$ is injective.
Thus we obtain an exact sequence $0 \to M \overset{f}{\to} C^l$ such that $(C^l)^{\dag} \overset{f^{\dag}}{\to} M^{\dag} \to 0$ is also exact.
Conversely, if there is an exact sequence $0 \to M \to C_0$ with $C_0\in\add\,C$ such that the $C$-dual sequence $C_0^{\dag} \to M^{\dag} \to 0$ is also exact, then dualizing the latter sequence by $C$, we get a commutative diagram
$$
\begin{CD}
0 @>>> M @>>> C_0 \\
@. @V{\lambda _M}VV @V{\lambda _{C_0}}VV \\
0 @>>> M^{\dag\dag} @>>> C_0^{\dag\dag}
\end{CD}
$$
with exact rows.
Noting that $C_0$ is a direct summand of $C^l=(R^l)^\dag$ for some $l>0$ and that $\lambda _{(R^l)^\dag}$ is injective by Remark \ref{triple}, we see that $\lambda _{C_0}$ is also injective.
Therefore $\lambda _M$ is also injective, that is, $M$ is $1$-$C$-torsionfree.

(2) According to (1), to show the assertion, we may assume that there is an exact sequence $0 \to M \to C_0 \to N \to 0$ with $C_0\in\add\,C$ whose $C$-dual sequence $0 \to N^\dag \to C_0^\dag \to M^\dag \to 0$ is also exact.
Dualizing the latter sequence by $C$ and using Remark \ref{triple}, we have a commutative diagram
$$
\begin{CD}
@. 0 @. 0 \\
@. @VVV @VVV \\
0 @>>> M @>>> C_0 @>>> N @>>> 0 \\
@. @V{\lambda _M}VV @V{\lambda _{C_0}}VV @V{\lambda _N}VV \\
0 @>>> M^{\dag\dag} @>>> C_0^{\dag\dag} @>>> N^{\dag\dag}
\end{CD}
$$
with exact rows and columns.
It follows from this diagram that if $M$ is $2$-$C$-torsionfree then $N$ is $1$-$C$-torsionfree, and that the converse holds when $\lambda _C$ is an isomorphism.
By (1), $N$ is $1$-$C$-torsionfree if and only if $N$ has an injective left $\add\,C$-approximation $N \to C_1$.
If this is the case, then the spliced exact sequence $0 \to M \to C_0 \to C_1$ is a left $\add\,C$-resolution.

(3) As we did in the proof of (2), we may assume that there is an exact sequence $0 \to M \to C_0 \to N \to 0$ with $C_0\in\add\,C$ such that the $C$-dual sequence $0 \to N^\dag \to C_0^\dag \to M^\dag \to 0$ is also exact.
We obtain a commutative diagram
$$
\begin{CD}
0 @>>> M @>>> C_0 @>>> N @>>> 0 \\
@. @V{\lambda _M}VV @V{\lambda _{C_0}}VV @V{\lambda _N}VV \\
0 @>>> M^{\dag\dag} @>>> C_0^{\dag\dag} @>>> N^{\dag\dag} \\
@>>> \Ext _R^1(M^\dag,C) @>>> \Ext _R^1(C_0^\dag,C) @>>> \Ext _R^1(N^\dag,C) \\
@>>> \cdots \\
@>>> \Ext _R^{n-2}(M^\dag, C) @>>> \Ext _R^{n-2}(C_0^\dag,C)
\end{CD}
$$
with exact rows.
From this diagram and the assumptions we see that if $M$ is $n$-$C$-torsionfree then $N$ is $(n-1)$-$C$-torsionfree, and that the converse holds when $\Ext _R^{n-2}(C^\dag,C)=0$.
By induction on $n$, $N$ is $(n-1)$-$C$-torsionfree if and only if it has an exact left $\add\,C$-resolution $0 \to N \to C_1 \to C_2 \to \cdots \to C_{n-1}$.
If this is the case, then the spliced exact sequence $0 \to M \to C_0 \to C_1 \to C_2 \to \cdots \to C_{n-1}$ is a left $\add\,C$-resolution of $M$.
\qed
\end{pf}

All the assumptions on $C$ in Proposition \ref{resol} are satisfied when $\Hom _R(C,C)$ is isomorphic to $R$:

\begin{cor}
Suppose that $\lambda _R$ is an isomorphism.
The following are equivalent for an $R$-module $M$.
\begin{enumerate}
\item[{\rm (1)}]
$M$ is $n$-$C$-torsionfree.
\item[{\rm (2)}]
$M$ has an exact left $\add\,C$-resolution $0 \to M \to C_0 \to C_1 \to \cdots \to C_{n-1}$.
\end{enumerate}
\end{cor}

\begin{pf}
Since $C^\dag\cong R$, the map $\lambda _C$ is an isomorphism and $\Ext _R ^i(C^\dag, C)=0$ for any $i>0$.
The assertion follows from Proposition \ref{resol}.
\qed
\end{pf}

We give here an application of Proposition \ref{resol}(1), which will be used later as a lemma.

\begin{cor}\label{1ctf}
Suppose that $\Ext _R ^1(C,C)=0$.
An $R$-module $M$ is $1$-$C$-torsionfree if and only if there is an exact sequence $0 \to M \to C_0 \to N \to 0$ such that $C_0 \in\add\, C$ and $\Ext _R ^1(N,C)=0$.
\end{cor}

\begin{pf}
Suppose that there is an exact sequence $0 \to M \to C_0 \to N \to 0$ such that $C_0\in\add\,C$ and $\Ext _R ^1(N,C)=0$.
Dualizing this sequence by $C$ gives an exact sequence $C_0^{\dag} \to M^{\dag} \to \Ext _R ^1 (N,C)=0$.
It follows from Proposition \ref{resol}(1) that $M$ is $1$-$C$-torsionfree.
Conversely, if this is the case, then we have an injective left $\add\,C$-approximation $f: M \to C_0$ by Proposition \ref{resol}(1).
Setting $N=\Coker\, f$, we get an exact sequence $0 \to M \overset{f}{\to} C_0 \to N \to 0$ such that $0 \to N^{\dag} \to C_0^{\dag} \overset{f^{\dag}}{\to} M^{\dag} \to 0$ is an exact sequence.
Since $\Ext _R ^1(C_0, C)=0$, we have $\Ext _R ^1 (N,C)=0$.
\qed
\end{pf}

Now, we can state and prove the following theorem, which is the main result of this paper.

\begin{thm}\label{main}
Let $C$ be an $n$-semidualizing $R$-module.
The following are equivalent for an $R$-module $M$.
\begin{enumerate}
\item[{\rm (1)}]
$\Omega ^n M$ is $n$-$C$-torsionfree.
\item[{\rm (2)}]
$M$ admits an $n$-$C$-spherical approximation.
\end{enumerate}
\end{thm}

\begin{pf}
Let $P_{\bullet}$ be a projective resolution of $M$.

(1) $\Rightarrow$ (2): We have an exact sequence $0 \to \Omega ^{i+1}M \to P_i \to \Omega ^iM \to 0$ for each $i$.
Set $X_0=\Omega ^nM$.
Note that $X_0$ is $n$-$C$-torsionfree.
Corollary \ref{1ctf} implies that there exists an exact sequence $0 \to X_0 \to C_0 \to Z_1 \to 0$ such that $C_0\in\add\,C$ and $\Ext _R ^1(Z_1, C)=0$.
We construct the pushout diagram:
$$
\begin{CD}
@. 0 @. 0 \\
@. @VVV @VVV \\
0 @>>> X_0 @>>> C_0 @>>> Z_1 @>>> 0 \\
@. @VVV @VVV @| \\
0 @>>> P_{n-1} @>>> X_1 @>>> Z_1 @>>> 0 \\
@. @VVV @VVV \\
@. \Omega ^{n-1}M @= \Omega ^{n-1}M \\
@. @VVV @VVV \\
@. 0 @. 0
\end{CD}
$$
Since $\Ext _R ^1(Z_1,C)=0=\Ext _R ^1(P_{n-1},C)$, we have $\Ext _R ^1(X_1,C)=0$.
If $n=1$, then the middle column is a desired exact sequence.

Let $n\geq 2$.
We establish the following two claims:

\setcounter{claim}{0}
\begin{claim}\label{Z_1}
$Z_1$ is $(n-1)$-$C$-torsionfree.
\end{claim}

\begin{cpf}
We have an exact sequence $0 \to Z_1^{\dag} \to C_0^{\dag} \to X_0^{\dag} \to 0$.
Dualize this again, and we have a commutative diagram
$$
\begin{CD}
0 @>>> X_0 @>>> C_0 @>>> Z_1 @>>> 0 \\
@. @V{\lambda _{X_0}}V{\cong}V @V{\lambda _{C_0}}VV @V{\lambda _{Z_1}}VV \\
0 @>>> X_0^{\dag\dag} @>>> C_0^{\dag\dag} @>{\eta}>> Z_1^{\dag\dag}
\end{CD}
$$
with exact rows.
Noting that $\lambda _R :R\to R^{\dag\dag}$ is an isomorphism, we see that so is $\lambda _C$, hence so is $\lambda _{C_0}$.
Thus the snake lemma says that $\lambda _{Z_1}$ is an injection.
Therefore, $Z_1$ is $(n-1)$-$C$-torsionfree when $n=2$.
When $n\geq 3$, noting that the map $\eta$ in the above diagram is a surjection because $\Ext _R ^1(X_0^{\dag}, C)=0$, we see from the five lemma that $\lambda _{Z_1}$ is an isomorphism.
Since $C_0$ belongs to $\add\,C$ and $\lambda _R: R\to C^\dag$ is an isomorphism, we easily see that the $R$-module $C_0^\dag$ is projective, and hence $\Ext _R ^i(C_0^\dag, C)=0$ for $i>0$.
Since $\Ext _R ^i(X_0^\dag, C)=0$ for $1\leq i\leq n-2$, we get $\Ext _R ^i(Z_1^\dag, C)=0$ for $1\leq i\leq n-3$.
It follows that $Z_1$ is $(n-1)$-$C$-torsionfree.
\qed
\end{cpf}

\begin{claim}\label{X_1}
$X_1$ is $(n-1)$-$C$-torsionfree.
\end{claim}

\begin{cpf}
We have an exact sequence $0 \to Z_1^{\dag} \to X_1^{\dag} \to P_{n-1}^{\dag} \to 0$, and dualizing this again yields a commutative diagram
$$
\begin{CD}
0 @>>> P_{n-1} @>>> X_1 @>>> Z_1 @>>> 0\\
@. @V{\lambda _{P_{n-1}}}V{\cong}V @V{\lambda _{X_1}}VV @V{\lambda _{Z_1}}VV \\
0 @>>> P_{n-1}^{\dag\dag} @>>> X_1^{\dag\dag} @>>> Z_1^{\dag\dag} @>>> 0
\end{CD}
$$
with exact rows.
If $n=2$, then $\lambda _{Z_1}$ is a monomorphism by Claim \ref{Z_1}, and so is $\lambda _{X_1}$.
Hence $X_1$ is $(n-1)$-$C$-torsionfree.
If $n\geq 3$, then $\lambda _{Z_1}$ is an isomorphism by Claim \ref{Z_1}, and so is $\lambda _{X_1}$.
Since $\Ext _R ^i(P_{n-1}^\dag, C)=0$ for $1\leq i\leq n$ and $\Ext _R ^i(Z_1^\dag, C)=0$ for $1\leq i\leq n-3$ by Claim \ref{Z_1}, we obtain $\Ext _R ^i(X_1^\dag, C)=0$ for $1\leq i\leq n-3$.
Thus $X_1$ is $(n-1)$-$C$-torsionfree.
\qed
\end{cpf}

According to Corollary \ref{1ctf}, there is an exact sequence $0 \to X_1 \to C_1 \to Z_2 \to 0$ with $C_1\in\add\,C$ and $\Ext _R ^1 (Z_2, C)=0$.
We construct the pushout diagram:
$$
\begin{CD}
@. 0 @. 0 \\
@. @VVV @VVV \\
@. C_0 @= C_0 \\
@. @VVV @VVV \\
0 @>>> X_1 @>>> C_1 @>>> Z_2 @>>> 0 \\
@. @VVV @VVV @| \\
0 @>>> \Omega ^{n-1}M @>>> Y_2 @>>> Z_2 @>>> 0 \\
@. @VVV @VVV \\
@. 0 @. 0
\end{CD}
$$
Using the bottom row of the above diagram, we construct the pushout diagram:
$$
\begin{CD}
@. 0 @. 0 \\
@. @VVV @VVV \\
0 @>>> \Omega ^{n-1}M @>>> Y_2 @>>> Z_2 @>>> 0 \\
@. @VVV @VVV @| \\
0 @>>> P_{n-2} @>>> X_2 @>>> Z_2 @>>> 0 \\
@. @VVV @VVV \\
@. \Omega ^{n-2}M @= \Omega ^{n-2}M \\
@. @VVV @VVV \\
@. 0 @. 0
\end{CD}
$$
From the first diagram, we immediately get $\Cdim\,Y_2<2$, and $\Ext _R ^2 (Z_2,C)=0$ because $\Ext _R ^1(X_1, C)=0=\Ext _R ^2(C_1, C)$.
Hence $\Ext _R ^i(Z_2, C)=0$ for $i=1, 2$, and we see from the middle row of the second diagram that $\Ext _R ^i(X_2,C)=0$ for $i=1,2$.
Thus, if $n=2$, then the middle column of the second diagram is a desired exact sequence.

Let $n\geq 3$.
Then similar arguments to the above claims show that both $Z_2$ and $X_2$ are $(n-2)$-$C$-torsionfree, and Corollary \ref{1ctf} yields an exact sequence $0 \to X_2 \to C_2 \to Z_3 \to 0$ such that $\Ext _R ^1(Z_3, C)=0$.
Similarly to the above, we construct two pushout diagrams:
$$
\begin{CD}
@. 0 @. 0 \\
@. @VVV @VVV \\
@. Y_2 @= Y_2 \\
@. @VVV @VVV \\
0 @>>> X_2 @>>> C_2 @>>> Z_3 @>>> 0 \\
@. @VVV @VVV @| \\
0 @>>> \Omega ^{n-2}M @>>> Y_3 @>>> Z_3 @>>> 0 \\
@. @VVV @VVV \\
@. 0 @. 0
\end{CD}
$$
and
$$
\begin{CD}
@. 0 @. 0 \\
@. @VVV @VVV \\
0 @>>> \Omega ^{n-2}M @>>> Y_3 @>>> Z_3 @>>> 0 \\
@. @VVV @VVV @| \\
0 @>>> P_{n-3} @>>> X_3 @>>> Z_3 @>>> 0 \\
@. @VVV @VVV \\
@. \Omega ^{n-3}M @= \Omega ^{n-3}M \\
@. @VVV @VVV \\
@. 0 @. 0
\end{CD}
$$
If $n=3$, then the middle column of the second diagram is a desired exact sequence.
If $n\geq 4$, then iterating this procedure, we eventually obtain an exact sequence $0 \to Y_n \to X_n \to M \to 0$ such that $\Ext _R ^i(X_n, C)=0$ for $1\leq i\leq n$ and $\Cdim\,Y_n<n$.

(2) $\Rightarrow$ (1): Let $0 \to Y \to X \to M \to 0$ be an $n$-$C$-spherical approximation of $M$.
Since $\Cdim\,Y<n$, there exists an exact sequence $0 \to C_{n-1} \overset{d_{n-1}}{\to} C_{n-2} \overset{d_{n-2}}{\to} \cdots \overset{d_1}{\to} C_0 \overset{d_0}{\to} Y \to 0$.
Put $Y_i=\Im\,d_i$ for each $i$.
We have exact sequences $0 \to Y_{i+1} \to C_i \to Y_i \to 0$ and $0 \to \Omega ^{i+1}M \to P_i \to \Omega ^iM \to 0$, where $P_i$ is projective $R$-module.
The following pullback diagram is obtained:
$$
\begin{CD}
@. @. 0 @. 0 \\
@. @. @VVV @VVV \\
@. @. \Omega M @= \Omega M \\
@. @. @VVV @VVV \\
0 @>>> Y @>>> L @>>> P_0 @>>> 0 \\
@. @| @VVV @VVV \\
0 @>>> Y @>>> X @>>> M @>>> 0 \\
@. @. @VVV @VVV \\
@. @. 0 @. 0
\end{CD}
$$
The projectivity of $P_0$ shows that the middle row splits; we have an isomorphism $L\cong Y\oplus P_0$.
Adding $P_0$ to the exact sequence $0 \to Y_1 \to C_0 \to Y \to 0$, we get an exact sequence $0 \to Y_1 \to C_0\oplus P_0 \to Y\oplus P_0 \to 0$.
Thus the following pullback diagram is obtained:
$$
\begin{CD}
@. 0 @. 0 \\
@. @VVV @VVV \\
@. Y_1 @= Y_1 \\
@. @VVV @VVV \\
0 @>>> X_1 @>>> C_0\oplus P_0 @>>> X @>>> 0 \\
@. @VVV @VVV @| \\
0 @>>> \Omega M @>>> Y\oplus P_0 @>>> X @>>> 0 \\
@. @VVV @VVV \\
@. 0 @. 0
\end{CD}
$$
Applying a similar argument to the left column of the above diagram, we get exact sequences $0 \to X_{i+1} \to C_i\oplus P_i \to X_i \to 0$ for $0\leq i\leq n-1$, where $X_0=X$ and $X_n=\Omega ^nM$.
The assumption yields $\Ext _R ^i(X_0,C)=0=\Ext _R ^i(C_0\oplus P_0, C)$ for $1\leq i\leq n$, hence we have an exact sequence $0 \to X_0^\dag \to (C_0\oplus P_0)^\dag \to X_1^\dag \to 0$ and $\Ext _R ^1(X_1, C)=0$ for $1\leq i\leq n-1$.
Inductively, for each $0\leq i\leq n-1$ an exact sequence $0 \to X_i^\dag \to (C_i\oplus P_i)^{\dag} \to X_{i+1}^\dag \to 0$ is obtained and $\Ext _R ^j (X_i, C)=0$ for $1\leq j\leq n-i$.
We have a commutative diagram
$$
\begin{CD}
0 @>>> X_1 @>>> C_0\oplus P_0 \\
@. @V{\lambda _{X_1}}VV @V{\lambda _{C_0\oplus P_0}}VV \\
0 @>>> X_1^{\dag\dag} @>>> (C_0\oplus P_0)^{\dag\dag}
\end{CD}
$$
with exact rows.
The assumption says that $\lambda _R$ is a monomorphism, and we see from Remark \ref{triple} that so is $\lambda _C=\lambda _{R^\dag}$.
Hence so is the map $\lambda _{C_0\oplus P_0}$, and so is $\lambda _{X_1}$.
Therefore $X_1$ is $1$-$C$-torsionfree.
If $n\geq 2$, then $\lambda _R$ is an isomorphism, and so is $\lambda _C$.
There is a commutative diagram
$$
\begin{CD}
0 @>>> X_2 @>>> C_1\oplus P_1 @>>> X_1 @>>> 0 \\
@. @V{\lambda _{X_2}}VV @V{\lambda _{C_1\oplus P_1}}VV @V{\lambda _{X_1}}VV \\
0 @>>> X_2^{\dag\dag} @>>> (C_1\oplus P_1)^{\dag\dag} @>>> X_1^{\dag\dag}
\end{CD}
$$
with exact rows.
Since $\lambda _{C_1\oplus P_1}$ is an isomorphism and $\lambda _{X_1}$ is a monomorphism, $\lambda _{X_2}$ is an isomorphism by the snake lemma.
Hence $X_2$ is $2$-$C$-torsionfree.
If $n\geq 3$, then we have a commutative diagram
$$
\begin{CD}
0 @>>> X_3 @>>> C_2\oplus P_2 @>>> X_2 @>>> 0 \\
@. @V{\lambda _{X_3}}VV @V{\lambda _{C_2\oplus P_2}}V{\cong}V @V{\lambda _{X_2}}V{\cong}V \\
0 @>>> X_3^{\dag\dag} @>>> (C_2\oplus P_2)^{\dag\dag} @>>> X_2^{\dag\dag} @>>> \Ext _R ^1(X_3^\dag, C) @>>> 0
\end{CD}
$$
with exact rows.
From this diagram it follows that $\lambda _{X_3}$ is an isomorphism and $\Ext _R ^1(X_3^\dag, C)=0$, which means that $X_3$ is $3$-$C$-torsionfree.
Repeating a similar argument, we see that $X_i$ is $i$-$C$-torsionfree for every $1\leq i\leq n$.
Therefore $\Omega ^nM=X_n$ is $n$-$C$-torsionfree, and the proof of the theorem is completed.
\qed
\end{pf}

\begin{rem}
The proof of Theorem \ref{main} actually shows a little stronger statements.
Let $M$ be an $R$-module.\\
(1) Assume that $\Ext _R ^1(C,C)=0$.
If $\Omega M$ is $1$-$C$-torsionfree, then $M$ has a $1$-$C$-spherical approximation.\\
(2) Assume that $\lambda _R$ is injective and $\Ext _R ^i(C, C)=0$ for $i=1, 2$.
If $\Omega ^2M$ is $2$-$C$-torsionfree, then $M$ has a $2$-$C$-spherical approximation.\\
(3) Assume that $\lambda _R$ is injective.
If $M$ has a $1$-$C$-spherical approximation, then $\Omega M$ is $1$-$C$-torsionfree.
(In fact, the injectivity of $\lambda _R$ implies the conclusion that $\Omega M$ is $1$-$C$-torsionfree by itself; see Lemma \ref{lambdasplit}(1) below.)\\
(4) Assume that $\lambda _R$ is bijective and $\Ext _R ^i(C, C)=0$ for $1\leq i\leq n-1$.
If $M$ has an $n$-$C$-spherical approximation, then $\Omega ^nM$ is $n$-$C$-torsionfree.
\end{rem}

\begin{rem}\label{aty}
Let $C$ be a semidualizing $R$-module (for the definition, see the following part of Definition \ref{sd}).
Recall the following:\\
(1) We say that an $R$-module $M$ is {\it totally $C$-reflexive} if $\lambda _M:M \to M^{\dag\dag}$ is an isomorphism and $\Ext _R ^i(M,C)=\Ext _R ^i(M^\dag, C)=0$ for any $i>0$.\\
(2) The {\it G$_C$-dimension} of an $R$-module $M$, which is denoted by $\GCdim _R\,M$, is defined as the infimum of integers $r$ such that there exists an exact sequence $0 \to X_r \to X_{r-1} \to \cdots \to X_0 \to M \to 0$ where each $X_i$ is totally $C$-reflexive.

Araya, Takahashi and Yoshino \cite{ATY} proved that every $R$-module of finite G$_C$-dimension admits an exact sequence which should be called an ``$\infty$-$C$-spherical approximation'' in our terminology; they proved that for every $R$-module $M$ with $\GCdim\,M<\infty$ there exists an exact sequence $0 \to Y \to X \to M \to 0$ of $R$-modules such that $X$ is totally $C$-reflexive and $\Cdim\,Y<\infty$.
\end{rem}

As a direct corollary of Theorem \ref{main}, we have an approximation theorem of Auslander and Bridger:

\begin{cor}
Let $M$ be an $R$-module.
Then $\Omega ^n M$ is $n$-torsionfree if and only if $M$ admits an $n$-spherical approximation.
\end{cor}

Theorem \ref{main} also induces a celebrated approximation theorem of Auslander and Buchweitz, which asserts the existence of Cohen-Macaulay approximations:

\begin{cor}
Let $(R, \m )$ be a $d$-dimensional Cohen-Macaulay local ring with the canonical module $\omega$.
Then every $R$-module $M$ admits a Cohen-Macaulay approximation, namely, there exists an exact sequence $0 \to Y \to X \to M \to 0$ such that $X$ is a maximal Cohen-Macaulay $R$-module and $Y$ is an $R$-module of finite injective dimension.
\end{cor}

\begin{pf}
If $d=0$, then $0 \to 0 \to M \overset{=}{\to} M \to 0$ is a desired exact sequence.
Let $d\geq 1$.
Then $\omega$ is $d$-semidualizing, and $\Omega ^dM$ is $d$-$\omega$-torsionfree.
Hence Theorem \ref{main} guarantees the existence of an exact sequence $0 \to Y \to X \to M \to 0$ such that $X$ is $d$-$\omega$-spherical and $\odim\,Y<d$.
Therefore $X$ is maximal Cohen-Macaulay.
On the other hand, noting that $\omega$ is an indecomposable $R$-module, we have an exact sequence $0 \to \omega ^{l_{d-1}} \to \omega ^{l_{d-2}} \to \cdots \to \omega ^{l_0} \to Y \to 0$.
Decomposing this into short exact sequences and noting that $\omega$ has finite injective dimension, one sees that $Y$ also has finite injective dimension.
\qed
\end{pf}

\section{Modules admitting approximations}

In the previous section, we proved that if $C$ is $n$-semidualizing, then any module whose $n^{\rm th}$ syzygy is $n$-$C$-torsionfree admits an $n$-$C$-spherical approximation.
In this section, we will consider modules whose $n^{\rm th}$ syzygies are $n$-$C$-torsionfree, and give several sufficient conditions for a given module to be such a module.

We begin by proving the following lemma.

\begin{lem}\label{lambdasplit}
Let $M$ be an $R$-module.
\begin{enumerate}
\item[{\rm (1)}]
If $R$ is $1$-$C$-torsionfree, then so is $\Omega M$.
\item[{\rm (2)}]
If $R$ is $2$-$C$-torsionfree, then for each $n\geq 2$ the map $\lambda _{\Omega ^n M}$ is a split monomorphism and the cokernel of $\lambda _{\Omega ^n M}$ is isomorphic to $\Ext _R ^n (M, C)^{\dag}$.
\end{enumerate}
\end{lem}

\begin{pf}
(1) There is an exact sequence $0 \to \Omega M \overset{\theta}{\to} P \to M \to 0$ with $P$ projective.
One has the following commutative diagram:
$$
\begin{CD}
\Omega M @>{\theta}>> P \\
@V{\lambda _{\Omega M}}VV @V{\lambda _{P}}VV \\
(\Omega M)^{\dag\dag} @>{\theta ^{\dag\dag}}>> P^{\dag\dag}
\end{CD}
$$
The map $\theta$ is injective, and the assumption implies that $\lambda _{P}$ is injective.
Hence $\lambda _{\Omega M}$ is also injective.

(2) Putting $X=\Omega ^{n-2}M$, one has $\Omega ^nM=\Omega ^2X$ and $\Ext _R ^n(M,C)=\Ext _R ^2(X,C)$.
Hence, replacing $M$ with $X$, we have only to show the lemma when $n=2$.
There is an exact sequence $0 \to \Omega ^2M \to Q \to \Omega M \to 0$, where $Q$ is a projective $R$-module.
Applying the $C$-dual functor $(-)^\dag$ yields an exact sequence $0 \to (\Omega M)^\dag \to Q^\dag \overset{\sigma}{\to} (\Omega ^2M)^\dag \to \Ext _R ^2(M,C) \to 0$.
Set $N=\Im\,\sigma$ and decompose this sequence into two short exact sequences:
\begin{equation}\label{2ses}
\begin{cases}
0 \to (\Omega M)^\dag \to Q^\dag \to N \to 0,\\
0 \to N \overset{\tau}{\to} (\Omega ^2M)^\dag \to \Ext _R ^2(M,C) \to 0.
\end{cases}
\end{equation}
Applying $(-)^\dag$ to the first sequence, one obtains a commutative diagram
$$
\begin{CD}
0 @>>> \Omega ^2M @>>> Q @>>> \Omega M @>>> 0 \\
@. @V{\eta}VV @V{\lambda _{Q}}VV @V{\lambda _{\Omega M}}VV \\
0 @>>> N^\dag @>>> Q^{\dag\dag} @>>> (\Omega M)^{\dag\dag}
\end{CD}
$$
with exact rows.
The map $\lambda _{Q}$ is an isomorphism by the assumption, and $\lambda _{\Omega M}$ is a monomorphism by (1).
Therefore we see from the snake lemma that the map $\eta$ is an isomorphism.
The diagram
$$
\begin{CD}
@. @. \Omega ^2M @= \Omega ^2M \\
@. @. @V{\lambda _{\Omega ^2M}}VV @V{\eta}V{\cong}V \\
0 @>>> \Ext _R ^2(M,C)^\dag @>>> (\Omega ^2M)^{\dag\dag} @>{\tau ^\dag}>> N^\dag
\end{CD}
$$
commutes, where the bottom row is the $C$-dual of the second sequence in \eqref{2ses}.
This commutative diagram shows that $\lambda _{\Omega ^2M}$ is a split monomorphism, that the sequence
\begin{equation}\label{later}
0 \to \Ext _R ^2(M,C)^\dag \to (\Omega ^2M)^{\dag\dag} \overset{\tau ^\dag}{\to} N^\dag \to 0
\end{equation}
is split exact and that $\Coker (\lambda _{\Omega ^2M})\cong\Ext _R ^2(M,C)^\dag$.
\qed
\end{pf}

For $R$-modules $M, N$, we define $\grade (M, N)$ by the infimum of integers $i$ such that $\Ext _R ^i(M, N)\neq 0$.
We state a criterion for $\Omega ^iM$ to be $i$-$C$-torsionfree for $1\leq i\leq n$ in terms of grade.

\begin{prop}\label{tfgr}
Let $C$ be an $R$-module such that $R$ is $(n-1)$-$C$-torsionfree.
\begin{enumerate}
\item[{\rm (1)}]
If $\Omega ^i M$ is $i$-$C$-torsionfree for every $1\leq i\leq n$, then $\grade (\Ext _R ^i(M, C), C)\geq i-1$ for every $1\leq i\leq n$.
\item[{\rm (2)}]
The converse also holds if $R$ is $n$-$C$-torsionfree.
\end{enumerate}
\end{prop}

\begin{pf}
We prove the proposition by induction on $n$.
When $n=1$, the statement (1) is trivial and the statement (2) follows from Lemma \ref{lambdasplit}(1).

Let $n=2$.
There is an exact sequence $0 \to \Omega ^2M \to P \to \Omega M \to 0$, where $P$ is projective.
Applying the functor $(-)^\dag$ gives another exact sequence $0 \to (\Omega M)^\dag \to P^\dag \to (\Omega ^2M)^\dag \to \Ext _R ^2(M,C) \to 0$.
Applying $(-)^\dag$ again, we get the following commutative diagram with exact rows:
$$
\begin{CD}
@. 0 @>>> \Omega ^2M @>>> P \\
@. @. @V{\lambda _{\Omega ^2M}}VV @V{\lambda _{P}}VV \\
0 @>>> \Ext _R ^2(M,C)^\dag @>>> (\Omega ^2M)^{\dag\dag} @>>> P^{\dag\dag}
\end{CD}
$$
Since $R$ is $1$-$C$-torsionfree by the assumption, $\lambda _R$ is injective, hence so is $\lambda _{P}$.
Thus, if $\Omega ^2M$ is $2$-$C$-torsionfree, then $\lambda _{\Omega ^2M}$ is bijective and diagram chasing shows $\Ext _R^2(M,C)^\dag =0$.
Therefore the statement (1) holds.
As for the statement (2), suppose that $R$ is $2$-$C$-torsionfree and $\Ext _R ^2(M,C)^\dag =0$.
The induction hypothesis shows that $\Omega M$ is $1$-$C$-torsionfree, and Lemma \ref{lambdasplit}(2) implies that $\lambda _{\Omega ^2M}$ is bijective, i.e., $\Omega ^2M$ is $2$-$C$-torsionfree.
Therefore the statement (2) holds.

Let $n\ge 3$.
The induction hypothesis implies that we may assume that $\Omega ^iM$ is $i$-$C$-torsionfree and $\grade (\Ext _R^i(M,C),C)\ge i-1$ for $1\le i\le n-1$.
There is an exact sequence
$$
0 \to \Omega ^nM \to P \to \Omega ^{n-1}M \to 0
$$
with $P$ projective, and dualizing this by $C$ yields an exact sequence $0 \to (\Omega ^{n-1}M)^\dag \to P^\dag \to (\Omega ^nM)^\dag \to \Ext _R ^n(M,C) \to 0$.
Decompose this sequence into two short exact sequences as follows:
$$
\begin{cases}
0 \to (\Omega ^{n-1}M)^\dag \to P^\dag \to N \to 0, \\
0 \to N \to (\Omega ^nM)^\dag \to \Ext _R ^n(M,C) \to 0.
\end{cases}
$$
Applying $(-)^\dag$ to the first exact sequence and using the assumption that $R$ is $(n-1)$-$C$-torsionfree, we get $\Ext _R^i(N,C)=0$ for $1\le i\le n-3$ and a monomorphism $\Ext _R^{n-2}(N,C) \hookrightarrow \Ext _R^{n-2}(P^\dag, C)$.
Apply $(-)^\dag$ to the second exact sequence.
Note that the proof of the existence of \eqref{later} in the proof of Lemma \ref{lambdasplit} shows that the induced sequence $0 \to \Ext _R ^n(M,C)^\dag \to (\Omega ^nM)^{\dag\dag} \to N^\dag \to 0$ is exact.
Hence we obtain isomorphisms $\Ext _R^i(\Ext _R^n(M,C),C)\cong\Ext _R^i((\Omega ^nM)^\dag, C)$ for $1\le i\le n-3$ and an exact sequence $0 \to \Ext _R^{n-2}(\Ext _R^n(M,C),C)\to\Ext _R^{n-2}((\Omega ^nM)^\dag, C)\to\Ext _R^{n-2}(N,C)$.
Moreover, Lemma \ref{lambdasplit}(2) guarantees that $\Ext _R^n(M,C)^\dag =0$ if and only if $\lambda _{\Omega ^nM}$ is an isomorphism.
Now it is easy to observe that if $\Omega ^nM$ is $n$-$C$-torsionfree then $\grade (\Ext _R^n(M,C),C)\ge n-1$, and that the converse holds when $R$ is $n$-$C$-torsionfree.
This completes the proof of the proposition.
\qed
\end{pf}

Next we want to consider the difference between the condition that $\Omega ^nM$ is $n$-$C$-torsionfree and the condition that $\Omega ^iM$ is $i$-$C$-torsionfree for $1\leq i\leq n$.
For this, let us study modules of finite $\add\,C$-resolution dimension.

\begin{lem}\label{cdimdepth}
Let $(R, \m, k)$ be a local ring and $C$ an $R$-module such that $\lambda _R$ is an isomorphism.
If $M$ is an $R$-module with $\Cdim\,M<\infty$, then $\Cdim\,M\leq\depth\,C$.
\end{lem}

\begin{pf}
First of all, let us show that the $R$-module $C$ is indecomposable: suppose that there is a direct sum decomposition $C\cong X\oplus Y$ for some nonzero $R$-modules $X,Y$.
Since $\lambda _R$ is an isomorphism, one has isomorphisms $R\cong\Hom _R (C,C)\cong\Hom _R (X,X)\oplus\Hom _R (X,Y)\oplus\Hom _R (Y,X)\oplus\Hom _R (Y,Y)$.
Both of the modules $\Hom _R (X,X)$ and $\Hom _R (Y,Y)$ are nonzero, because they contain the identity maps.
Hence $R$ is decomposable as an $R$-module, which contradicts the assumption that $R$ is a local ring.
Therefore $C$ is an indecomposable $R$-module, and all objects of $\add\,C$ are finite direct sums of copies of $C$.

Put $\Cdim\,M=s$.
There exists an exact sequence $0 \to C^{l_s} \overset{\phi _s}{\to} C^{l_{s-1}} \overset{\phi _{s-1}}{\to} \cdots \overset{\phi _1}{\to} C^{l_0} \to M \to 0$.
Note that all $l_i$ are nonzero.
Since the homothety map $\lambda _R: R\to \Hom _R (C,C)$ is an isomorphism, one can identify each $\phi _i$ with an $l_{i-1}\times l_i$ matrix over $R$.
With respect to the suitable bases of $C^{l_s}$ and $C^{l_{s-1}}$, the matrix $\phi _s$ has the form $\left(
\begin{smallmatrix}
E & 0 \\
0 & A
\end{smallmatrix}
\right)$, where $E$ is an identity matrix and $A$ is a matrix whose components are in the maximal ideal $\m$.
Removing $E$ from $\phi _s$ and iterating this procedure, we may assume that all components of the matrix $\phi _i$ belong to the maximal ideal $\m$ for each $1\leq i\leq s$.

Set $t=\depth\,C$.
We want to prove that $s$ is not bigger than $t$.
If $s=0$, then the inequality obviously holds.
Let $s>0$.
From the exact sequence $0 \to C^{l_s} \overset{\phi _s}{\to} C^{l_{s-1}}$ we get an exact sequence $0 \to \Hom _R (k,C^{l_s}) \overset{f}{\to} \Hom _R (k,C^{l_{s-1}})$, where $f=\Hom _R (k,\phi _s)$.
Since all the components of $\phi _s$ belong to $\m$, we have $f=0$.
Hence $\Hom _R (k,C^{l_s})=0$, which implies that $t=\depth\,C>0$.
Putting $M_i=\Im\,\phi _i$, we have an exact sequence $0 \to M_{i+1} \to C^{l_i} \to M_i \to 0$ for $1\leq i\leq s-1$.
Assume that $s$ is bigger than $t$.
Then from the exact sequence $0 \to C^{l_s} \overset{\phi _s}{\to} C^{l_{s-1}} \to M_{s-1} \to 0$ we get an exact sequence $0=\Ext _R ^{t-1}(k,C^{l_{s-1}}) \to \Ext _R ^{t-1}(k, M_{s-1}) \to \Ext _R ^t(k, C^{l_s}) \overset{g}{\to} \Ext _R ^t(k, C^{l_{s-1}})$, where $g=\Ext _R ^t(k,\phi _s)$.
Since all the components of the matrix $\phi _s$ are in $\m$, the map $g$ is a zero map.
Noting that $\Ext _R ^t(k,C)\neq 0$, we have $\Ext _R ^{t-1}(k,M_{s-1})\neq 0$.
There are isomorphisms $\Ext _R ^{t-1}(k,M_{s-1})\cong\Ext _R ^{t-2}(k,M_{s-2})\cong\cdots\cong\Ext _R ^1(k,M_{s-t+1})\cong\Hom _R (k,M_{s-t})$, which show that $\depth\,M_{s-t}=0$.
However, there is an injection $M_{s-t}\to C^{l_{s-t-1}}$; note that $s-t-1\geq 0$.
Hence $t=\depth\,C^{l_{s-t-1}}=0$, which is a contradiction.
This contradiction proves that $s$ is not bigger than $t$.
\qed
\end{pf}

\begin{lem}\label{nmc}
Let $r$ be an integer, and let $C$ be an $R$-module such that $\Ext _R ^i(C,C)=0$ for all $1\leq i\leq r$.
If $M$ is an $R$-module with $\Cdim\,M<r$, then one has $\Ext _R ^r(M,C)=0$.
\end{lem}

\begin{pf}
Putting $s=\Cdim\,M$, we have an exact sequence $0 \to C_s \to C_{s-1} \to \cdots \to C_0 \to M \to 0$ such that each $C_i$ belongs to $\add\,C$.
Decomposing this into short exact sequences, we get an exact sequence $0 \to M_{i+1} \to C_i \to M_i \to 0$ for each $1\leq i\leq s-1$.
There are isomorphisms $\Ext _R ^r(M,C)\cong\Ext _R ^{r-1}(M_1,C)\cong\cdots\cong\Ext _R ^{r-s+1}(M_{s-1},C)\cong\Ext _R ^{r-s}(C_s,C)$, where the last Ext module vanishes because $1\leq r-s\leq r$.
\qed
\end{pf}

The following is a well-known result concerning grade; see \cite[Proposition 1.2.10(a),(e)]{BH}, for example.

\begin{lem}\label{ae}
For $R$-modules $M$ and $N$, one has $\grade (M,N)=\inf\{\,\depth\,N_\p\,|\,\p\in\Supp\,M\,\}$.
\end{lem}

Combining the above three lemmas, we obtain the following.

\begin{lem}\label{cdimgrade}
Let $C$ be an $R$-module such that $\lambda _R$ is an isomorphism and $\Ext _R ^i(C,C)=0$ for $1\leq i\leq n$.
If $M$ is an $R$-module with $\Cdim\,M<\infty$, then $\grade (\Ext _R ^i(M,C),C)\geq i$ for all $1\leq i\leq n$.
\end{lem}

\begin{pf}
Fix an integer $i$ with $1\leq i\leq n$.
Let $\p$ be an prime ideal of $R$ satisfying $\depth _{R_\p}\,C_\p <i$.
Then one has $\Cpdim _{R_\p}\,M_\p <i$ by Lemma \ref{cdimdepth}, and hence $\Ext _{R_\p}^i(M_\p , C_\p )=0$ by Lemma \ref{nmc}.
Therefore $\p\not\in\Supp\,\Ext _R^i(M,C)$, and the inequality $\grade (\Ext _R^i(M,C),C)\geq i$ follows from Lemma \ref{ae}.
\qed
\end{pf}

Now we can consider the relationship between the $n$-$C$-torsionfreeness of $\Omega ^nM$ and the $i$-$C$-torsionfreeness of $\Omega ^iM$ for $1\le i\le n$.
Under the assumption that $C$ is $n$-semidualizing, these properties are equivalent to each other.
We should remark that an $n$-$C$-spherical approximation of $M$ plays an essential role in the proof.

\begin{prop}\label{descent}
Let $C$ be an $n$-semidualizing $R$-module.
The following are equivalent for an $R$-module $M$.
\begin{enumerate}
\item[{\rm (1)}]
$\Omega ^nM$ is $n$-$C$-torsionfree.
\item[{\rm (2)}]
$\Omega ^iM$ is $i$-$C$-torsionfree for every $1\leq i\leq n$.
\end{enumerate}
\end{prop}

\begin{pf}
(2) $\Rightarrow$ (1): This implication is trivial.

(1) $\Rightarrow$ (2): Proposition \ref{tfgr} says that it is enough to prove that the inequality $\grade (\Ext _R ^i(M,C),C)\geq i-1$ holds for $1\leq i\leq n$.
This inequality automatically holds when $n=1$, hence let $n\geq 2$.
Theorem \ref{main} yields an $n$-$C$-spherical approximation $0 \to Y \to X \to M \to 0$ of $M$.
Hence we have an isomorphism $\Ext _R ^i(M,C)\cong\Ext _R ^{i-1}(Y,C)$ for $2\leq i\leq n$.
Since $n\geq 2$, the assumption that $C$ is $n$-semidualizing implies that $\lambda _R$ is an isomorphism and $\Ext _R ^i(C,C)=0$ for $1\leq i\leq n$.
Therefore we see from Lemma \ref{cdimgrade} that $\grade (\Ext _R ^i(Y,C),C)\geq i$ for $1\leq i\leq n$, and thus $\grade (\Ext _R ^i(M,C),C)\geq i-1$ for $1\leq i\leq n$, as desired.
\qed
\end{pf}

Our next aim is to give sufficient conditions for an $R$-module to satisfy the condition that the $n^{\rm th}$ syzygy is $n$-$C$-torsionfree.
For this aim, we need to introduce the following lemma which will also be used later, and to recall basic properties of G$_C$-dimension.
(For the definition of G$_C$-dimension, see Remark \ref{aty}.)

\begin{lem}\label{1441}
Let $R$ be a local ring and $r$ a positive integer.
Suppose that $\lambda _R$ is an isomorphism and $\Ext _R ^i(C,C)=0$ for all $1\leq i<r$.
Then the following hold.
\begin{enumerate}
\item[{\rm (1)}]
$\depth\,R\geq r$ if and only if $\depth\,C\geq r$.
\item[{\rm (2)}]
Let $R$ be a Cohen-Macaulay local ring with $\Kdim\,R<r$.
Then $C$ is a maximal Cohen-Macaulay $R$-module.
\end{enumerate}
\end{lem}

\begin{pf}
(1) Let us show the assertion by induction on $r$.
Since $R\cong\Hom _R(C,C)$, one has $\Hom _R(k,R)\cong\Hom _R(k,\Hom _R(C,C))\cong\Hom _R(k\otimes _R C, C)\cong\Hom _R(k,C)^m$, where $m$ is the minimal number of generators of $C$.
Hence $\Hom _R(k,R)=0$ if and only if $\Hom _R(k,C)=0$.
In other words, $\depth\,R\geq 1$ if and only if $\depth\,C\geq 1$.
Thus the assertion in the case where $r=1$ is proved.

Let $r\geq 2$.
The induction hypothesis guarantees the existence of an element $x\in R$ which is both $R$-regular and $C$-regular.
Set $\overline{(-)}=(-)\otimes _R R/(x)$.
There is an exact sequence $0 \to C \overset{x}{\to} C \to \overline C \to 0$.
Apply the functor $(-)^\dag$ to this, and make a long exact sequence; we see that the map $\overline R\to\Ext _R^1(\overline C,C)$ induced by the connecting homomorphism $R\cong C^\dag \to \Ext _R ^1 (\overline C,C)$ is an isomorphism, and that $\Ext _R ^i(\overline C, C)=0$ for $2\leq i<r$.
Since there is a natural isomorphism $\Ext _{\overline R} ^i(\overline C,\overline C)\cong\Ext _R ^{i+1}(\overline C,C)$ for any integer $i$, the natural map $\overline R \to \Hom _{\overline R}(\overline C, \overline C)$ is an isomorphism and $\Ext _{\overline R}^i(\overline C, \overline C)=0$ for $1\leq i<r-1$.
It follows from the induction hypothesis that $\depth\,\overline R\geq r-1$ if and only if $\depth\,\overline C\geq r-1$.
Thus $\depth\,R\geq r$ if and only if $\depth\,C\geq r$.

(2) Set $d=\Kdim\,R$.
According to the assumptions, the map $\lambda _R$ is an isomorphism, $\Ext _R ^i(C,C)=0$ for $1\leq i\leq d$ and $\depth\,R=d\geq d$.
Hence the assertion (1) yields $\depth\,C\geq d$, which means that $C$ is maximal Cohen-Macaulay.
\qed
\end{pf}

\begin{lem}\cite{Golod}\label{gcdim}
Let $R$ be a local ring and $C$ a semidualizing $R$-module.
Then the following hold for an $R$-module $M$:
\begin{enumerate}
\item[{\rm (1)}]
If $\GCdim _R\,M<\infty$, then $\GCdim _R\,M=\depth\,R-\depth _R\,C=\sup\{\,i\,|\,\Ext _R^i(M,C)\neq 0\,\}$,
\item[{\rm (2)}]
For a nonnegative integer $r$, one has $\GCdim _R (\Omega ^rM)=\sup\{\,\GCdim _R\,M-r,\ 0\,\}$.
\end{enumerate}
\end{lem}

Let us give sufficient conditions for an $R$-module $M$ to be such that $\Omega ^nM$ is $n$-$C$-torsionfree:

\begin{prop}\label{condm}
Let $M$ be an $R$-module, and let $C$ be an $R$-module such that $R$ is $n$-$C$-torsionfree.
Suppose that either of the following holds:
\begin{enumerate}
\item[{\rm (1)}]
$\pd _{R_{\p}}\, M_{\p}<\infty$ for any $\p\in\Spec\,R$ with $\depth\,R_{\p}\leq n-2$,
\item[{\rm (2)}]
$C_{\p}$ is a semidualizing $R_{\p}$-module and $\GCpdim _{R_{\p}}\, M_{\p}<\infty$ for any $\p\in\Spec\,R$ with $\depth\,R_{\p}\leq n-2$.
\end{enumerate}
Then $\Omega ^n M$ is $n$-$C$-torsionfree.
\end{prop}

\begin{pf}
The proposition can be proved by induction on $n$.
When $n=1$, by Lemma \ref{lambdasplit}(1), the conclusion automatically holds.
Hence let $n\geq 2$.
It follows from the induction hypothesis that $\Omega ^iM$ is $i$-$C$-torsionfree for $1\leq i\leq n-1$.
Proposition \ref{tfgr} says that $\grade (\Ext _R ^i(M,C),C)\geq i-1$ for $1\leq i\leq n-1$, and that we have only to prove the inequality
$$
\grade (\Ext _R ^n(M,C),C)\geq n-1.
$$
Let $\p$ be a prime ideal of $R$ satisfying $\depth\,C_\p\leq n-2$.
Then $\depth\,R_\p\leq n-2$ by Lemma \ref{1441}(1).

Firstly, suppose that the condition (1) holds.
Then $\pd _{R_\p}\,M_\p <\infty$, and one has $\pd _{R_\p}\,M_\p\leq\depth\,R_\p\leq n-2$.
Hence $\Ext _{R_\p}^i(M_\p, X)=0$ for any $R_\p$-module $X$ and any $i>n-2$.
In particular, we get $\Ext _R ^n(M,C)_\p\cong\Ext _{R_\p}^n(M_\p,C_\p )=0$, and Lemma \ref{ae} yields the inequality $\grade (\Ext _R ^n(M,C),C)\geq n-1$.

Secondly, suppose that the condition (2) holds.
Then $\GCpdim _{R_\p}\,M_\p<\infty$, and Lemma \ref{gcdim}(1) implies that $\sup\{\,i\,|\,\Ext _{R_\p}^i(M_\p,C_\p)\neq 0\,\}=\GCpdim _{R_\p}\,M_\p\leq\depth\,R_\p\leq n-2$.
Therefore one has $\Ext _{R_\p}^n(M_\p,C_\p)=0$.
It follows from Lemma \ref{ae} that the inequality $\grade (\Ext _R ^n(M,C),C)\geq n-1$ holds.
\qed
\end{pf}

Recall that the G-dimension of an $R$-module $M$, which is denoted by $\Gdim _R\,M$, is defined as the G$_R$-dimension of $M$.
As a corollary of the above proposition, we get a theorem of Ma\c{s}ek.

\begin{cor}\cite[Theorem 43]{Masek}
Let $X$ be an $R$-module such that $\Gdim _{R_\p}\,X_\p<\infty$ for any $\p\in\Spec\,R$ with $\depth\,R_\p\leq n-2$.
Then $X$ is $n$-torsionfree if (and only if) $X$ is $n$-syzygy.
\end{cor}

\begin{pf}
Suppose that $X$ is $n$-syzygy, i.e., $X\cong\Omega ^nM$ for some $R$-module $M$.
Then note by Lemma \ref{gcdim}(2) that the $R_\p$-module $M_\p$ is also of finite G-dimension for $\p\in\Spec\,R$ with $\depth\,R_\p\leq n-2$.
Proposition \ref{condm} shows that $\Omega ^nM$ is $n$-torsionfree.
Hence $X$ is $n$-torsionfree.
\qed
\end{pf}

\section{Contravariant finiteness}

In the previous section, we investigated those $R$-modules whose $n^{\rm th}$ syzygies are $n$-$C$-torsionfree; we obtained several conditions for a {\em given} $R$-module to be such a module.
In this section, we shall consider conditions for {\em all} $R$-modules to be such modules.
After that, we will study the contravariant finiteness of the full subcategory of $\mod\,R$ consisting of all $n$-$C$-spherical modules.

First of all, we prove the following theorem, which is the second main result of this paper.

\begin{thm}\label{cond}
Suppose that $R$ is $n$-$C$-torsionfree.
Then the following are equivalent.
\begin{enumerate}
\item[{\rm (1)}]
$\id _{R_{\p}}\, C_{\p}<\infty$ for any $\p\in\Spec\,R$ with $\depth\,R_{\p}\leq n-2$.
\item[{\rm (2)}]
$\Omega ^n M$ is $n$-$C$-torsionfree for any $R$-module $M$.
\end{enumerate}
\end{thm}

\begin{pf}
We use induction on $n$ to prove the theorem.
When $n=1$, the assertion (1) holds because there is no prime ideal $\p$ of $R$ satisfying $\depth\,R_{\p}\leq n-2$, and the assertion (2) holds by Lemma \ref{lambdasplit}(1).
In the following, we consider the case where $n\geq 2$.

(1) $\Rightarrow$ (2): Fix an $R$-module $M$.
The induction hypothesis shows that $\Omega ^iM$ is $i$-$C$-torsionfree for $1\leq i\leq n-1$.
By Proposition \ref{tfgr}, we have $\grade (\Ext _R ^i(M,C),C)\geq i-1$ for $1\leq i\leq n-1$, and it suffices to prove that the inequality $\grade (\Ext _R ^n(M,C),C)\geq n-1$ holds.
Let $\p\in\Spec\,R$.
If $\depth\,C_\p\leq n-2$, then $\depth\,R_\p\leq n-2$ by Lemma \ref{1441}(1).
The assumption says that $\id _{R_\p}\,C_\p<\infty$, and $\id _{R_\p}\,C_\p=\depth\,R_\p\leq n-2$.
Therefore $\Ext _R^n(M,C)_\p\cong\Ext _{R_\p}^n(M_\p, C_\p)=0$.
Thus we see from Lemma \ref{ae} that $\grade (\Ext _R ^n(M,C),C)\geq n-1$, as desired.

(2) $\Rightarrow$ (1): When $n=2$, Lemma \ref{lambdasplit}(2) implies that $\Ext _{R_\p}^2(M_\p , C_\p )=0$ for all $R$-modules $M$ and $\p\in\Ass\,C$, because $\Ass (\Ext _R ^2(M,C)^\dag )=\Supp\,\Ext _R ^2(M,C)\cap\Ass\,C$.
The isomorphism $\lambda _R:R\to \Hom _R (C,C)$ shows that $\Ass\,C$ coincides with $\Ass\,R$.
Hence, setting $M=\Omega _R ^i(R/\p)$, one has $\Ext _{R_\p} ^{i+2} (\kappa (\p) , C_\p )\cong\Ext _{R_\p} ^2((\Omega _R^i(R/\p ))_\p ,C_\p )=0$ for any $\p\in\Ass\,R$ and any $i>0$.
Therefore $\id _{R_\p}\,C_\p <\infty$ for $\p\in\Spec\,R$ with $\depth\,R_\p =0$.

Let $n\geq 3$.
Fix an $R$-module $M$.
We have an exact sequence $0 \to \Omega ^{n+1}M \to P \to \Omega ^nM \to 0$ such that $P$ is a projective $R$-module.
From this we get another exact sequence $0 \to (\Omega ^nM)^\dag \to P^\dag \to (\Omega ^{n+1}M)^\dag \to \Ext _R ^{n+1} (M,C) \to 0$.
Decompose this into short exact sequences:
\begin{equation}\label{twoses}
\begin{cases}
0 \to (\Omega ^nM)^\dag \to P^\dag \to N \to 0, \\
0 \to N \to (\Omega ^{n+1}M)^\dag \to \Ext _R ^{n+1} (M,C) \to 0.
\end{cases}
\end{equation}
Note from the assumption that both $\Omega ^nM$ and $\Omega ^{n+1}M=\Omega ^n(\Omega M)$ are $n$-$C$-torsionfree.
Since $R$ is $n$-$C$-torsionfree, we see from the first sequence in \eqref{twoses} that there is a commutative diagram
$$
\begin{CD}
0 @>>> \Omega ^{n+1}M @>>> P @>>> \Omega ^nM @>>> 0 \\
@. @V{\alpha}VV @V{\lambda _{P}}V{\cong}V @V{\lambda _{\Omega ^nM}}V{\cong}V \\
0 @>>> N^\dag @>>> P^{\dag\dag} @>>> (\Omega ^nM)^{\dag\dag} @>>> \Ext _R ^1(N,C) @>>> 0
\end{CD}
$$
with exact rows, and $\Ext _R ^i(N,C)=0$ for $2\leq i\leq n-2$.
This diagram shows that $\alpha$ is an isomorphism and $\Ext _R ^1(N,C)=0$.
The second sequence in \eqref{twoses} gives an exact sequence $0 \to \Ext _R ^{n+1} (M,C)^\dag \to (\Omega ^{n+1}M)^{\dag\dag} \overset{\beta}{\to} N^\dag \to \Ext _R ^1(\Ext _R ^{n+1}(M,C),C) \to 0$ and $\Ext _R ^i(\Ext _R ^{n+1}(M,C),C)=0$ for $2\leq i\leq n-2$.
Since the diagram
$$
\begin{CD}
\Omega ^{n+1} M @= \Omega ^{n+1}M \\
@V{\lambda _{\Omega ^{n+1}M}}V{\cong}V @V{\alpha}V{\cong}V \\
(\Omega ^{n+1}M)^{\dag\dag} @>{\beta}>> N^\dag
\end{CD}
$$
commutes, the map $\beta$ is an isomorphism, and $\Ext _R ^{n+1} (M,C)^\dag=0=\Ext _R ^1(\Ext _R ^{n+1}(M,C),C)$.
Thus we have $\Ext _R ^i(\Ext _R ^{n+1}(M,C),C)=0$ for every $i\leq n-2$, which means that the inequality $\grade (\Ext _R ^{n+1} (M,C), C) \geq n-1$ holds.
Therefore, if $\p$ is a prime ideal of $R$ with $\depth\,R_\p\leq n-2$, then $\depth\,C_\p \leq n-2$ by Lemma \ref{1441}(1), and it follows from Lemma \ref{ae} that $\p$ does not belong to $\Supp\,\Ext _R^{n+1}(M,C)$, i.e., $\Ext _{R_\p}^{n+1} (M_\p , C_\p )=0$.
Putting $M=\Omega _R ^i(R/\p)$, we obtain $\Ext _{R_\p}^{n+1+i} (\kappa (\p), C_\p )\cong\Ext _{R_\p}^{n+1}((\Omega _R^i(R/\p ))_\p, C_\p)=0$ for any $i>0$.
This implies that $\id _{R_\p}\,C_\p<\infty$, and the proof is completed.
\qed
\end{pf}

\begin{rem}
In the case where $C$ is $n$-semidualizing, one can prove the above theorem more easily, as follows.
Fix an $R$-module $M$.
The $n$-$C$-torsionfree property of $\Omega ^nM$ is equivalent to the $i$-$C$-torsionfree property of $\Omega ^iM$ for $1\leq i\leq n$ by Proposition \ref{descent}, which is equivalent to the inequality $\grade (\Ext _R ^i(M,C),C)\geq i-1$ for $1\leq i\leq n$ by Proposition \ref{tfgr}, which is equivalent to the condition that $\Ext _R ^i(M_\p, C_\p )=0$ for any $\p\in\Spec\,R$ with $\depth\,R_\p\leq i-2$ by Lemmas \ref{ae} and \ref{1441}(1).
Hence, setting $M=\Omega ^j(R/\p )$ for $j\gg 0$, we easily see that the two conditions in the theorem are equivalent.
\end{rem}

\begin{rem}\label{suff}
It is easy to see that Theorem \ref{cond}(1) holds in each of the following cases:\\
(1) $R$ satisfies Serre's condition $(S_{n-1})$ and $C$ is locally of finite injective dimension in codimension $n-2$, namely, $\id _{R_\p}\,C_\p <\infty$ for $\p\in\Spec\,R$ with $\height\,\p\leq n-2$.\\
(2) $(R, \m )$ is local, $n$ is at most $\depth\,R+1$ and $C$ is locally of finite injective dimension on the punctured spectrum of $R$, namely, $\id _{R_\p}\,C_\p <\infty$ for $\p\in\Spec\,R-\{\m\}$.
\end{rem}

\begin{cor}
Let $C$ be an $n$-semidualizing $R$-module.
Suppose that $R$ satisfies Serre's conditions $(R_{n-2})$ and $(S_{n-1})$.
Then every $R$-module $M$ has an $n$-$C$-spherical approximation.
\end{cor}

\begin{pf}
Take a prime ideal $\p$ of $R$ with $\height\,\p \leq n-2$.
Then $R_\p$ is a regular local ring, hence $\id _{R_\p}\,C_\p <\infty$.
Thus the corollary follows from Remark \ref{suff}(1) and Theorems \ref{cond} and \ref{main}.
\qed
\end{pf}

The lemma below says that over a Gorenstein local ring of dimension $d\geq 2$, any $n$-semidualizing module is free for $n\geq d$.

\begin{lem}\label{c=r}
Let $(R, \m , k)$ be a $d$-dimensional Gorenstein local ring.
If $\lambda _R$ is an isomorphism and $\Ext _R ^i(C,C)=0$ for $1\leq i\leq d$, then $C\cong R$.
\end{lem}

\begin{pf}
Lemma \ref{1441}(2) says that the $R$-module $C$ is maximal Cohen-Macaulay.
Hence there exists a sequence $\xx =x_1,\dots,x_d$ in $\m$ which is both $R$-regular and $C$-regular.
Repeating a similar argument to the proof of Lemma \ref{1441}(1), we obtain
$$
\Ext _{R/(\xx _i)} ^j(C/\xx _iC, C/\xx _iC)\cong
\begin{cases}
R/(\xx _i) & \text{if } j=0, \\
0 & \text{if } 1\leq j\leq d-i,
\end{cases}
$$
where $\xx _i=x_1,\dots,x_i$, for each $1\leq i\leq d$.
It follows that there is an isomorphism $\Hom _{\overline R} (\overline C, \overline C)\cong \overline R$, where $\overline{(-)}=(-)\otimes _R R/(\xx)$.
Hence $\Hom _{\overline R}(k,\overline R)\cong\Hom _{\overline R}(k, \Hom _{\overline R}(\overline C,\overline C))\cong\Hom _{\overline R}(k\otimes _{\overline R}\overline C,\overline C)\cong\Hom _{\overline R}(k, \overline C)^m$, where $m$ is the minimal number of generators of the $\overline R$-module $\overline C$.
Since $\overline R$ is artinian Gorenstein, we have $\Hom _{\overline R} (k,\overline R)\cong k$.
Thus one must have $m=1$, hence $\overline C$ is a cyclic $\overline R$-module, and $\overline C\cong \overline R/I$ for some ideal $I$ of $\overline R$.
Therefore $\overline R\cong\Hom _{\overline R}(\overline C,\overline C)\cong\Hom _{\overline R}(\overline R/I,\overline R/I)\cong\overline R/I$, and we get $I=0$, i.e., $\overline C\cong\overline R$.
This gives an isomorphism $C\cong R$.
\qed
\end{pf}

Applying the above lemma, we can get a sufficient condition for $R$ and $C$ to satisfy Theorem \ref{cond}(1).

\begin{prop}
Suppose that $R$ is $n$-$C$-torsionfree and that $R_\p$ is Gorenstein for any $\p\in\Spec\,R$ with $\depth\,R_{\p}\leq n-2$.
Then $\id _{R_{\p}}\, C_{\p}<\infty$ for any $\p\in\Spec\,R$ with $\depth\,R_{\p}\leq n-2$.
(Hence $\Omega ^nM$ is $n$-$C$-torsionfree for any $R$-module $M$.)
\end{prop}

\begin{pf}
We may assume $n\geq 2$.
Take a prime ideal $\p$ of $R$ with $\depth\,R_\p\leq n-2$.
Then $R_\p$ is Gorenstein, and $\Kdim\,R_\p = \depth\,R_\p\leq n-2$.
Since $R$ is $n$-$C$-torsionfree, one has
$$
\Ext _{R_\p}^i(C_\p, C_\p)\cong
\begin{cases}
R_\p & \text{if } i=0, \\
0 & \text{if }1\leq i\leq \Kdim\,R_\p.
\end{cases}
$$
Hence Lemma \ref{c=r} yields an isomorphism $C_\p\cong R_\p$.
Therefore $\id _{R_\p}\,C_\p=\id _{R_\p}\,R_\p<\infty$.
\qed
\end{pf}

Here we notice that an $n$-$C$-spherical approximation gives a right approximation:

\begin{prop}\label{xy}
Define two full subcategories of $\mod\,R$ as follows:
\begin{align*}
\X & = \{\,X\in\mod\,R\,|\,X\text{ is }n\text{-}C\text{-spherical}\,\},\\
\Y & =\{\,Y\in\mod\,R\,|\,\Cdim\,Y<n\,\}.
\end{align*}
Let $0 \to Y \to X \overset{f}{\to} M \to 0$ be an exact sequence of $R$-modules with $X\in\X$ and $Y\in\Y$.
Then the homomorphism $f$ is a right $\X$-approximation of $M$.
\end{prop}

\begin{pf}
One has $\Ext _R ^1(X',Y')=0$ for any $X'\in\X$ and $Y'\in\Y$.
In fact, there exists an exact sequence $0 \to C_{n-1} \overset{d_{n-1}}{\to} C_{n-2} \overset{d_{n-2}}{\to} \cdots \overset{d_1}{\to} C_0 \overset{d_0}{\to} Y' \to 0$.
Putting $Y_i=\Im\,d_i$, we have an exact sequence $0 \to Y_{i+1} \to C_i \to Y_i \to 0$ for each $i$.
Since $\Ext _R ^i(X',C)=0$ for $1\leq i\leq n$, we obtain isomorphisms $\Ext _R ^1(X',Y')\cong\Ext _R ^2(X',Y_1)\cong\Ext _R ^3(X',Y_2)\cong\cdots\cong\Ext _R ^{n-1}(X',Y_{n-2})\cong\Ext _R ^n(X',C_{n-1})=0$, as desired.
For any $X'\in\X$, we have an exact sequence $\Hom _R (X',X) \overset{g}{\to} \Hom _R (X',M) \to \Ext _R ^1(X',Y)=0$ where $g=\Hom _R (X',f)$.
This says that the homomorphism $f$ is a right $\X$-approximation of $M$.
\qed
\end{pf}

Using Theorems \ref{cond} and \ref{main} and Proposition \ref{xy}, we see that the subcategory of $n$-$C$-spherical $R$-modules is contravariantly finite under some assumptions:

\begin{cor}
Let $C$ be an $n$-semidualizing $R$-module such that $\id _{R_{\p}}\, C_{\p}<\infty$ for any $\p\in\Spec\,R$ with $\depth\,R_{\p}\leq n-2$.
Then the full subcategory
$$
\X = \{\,X\in\mod\,R\,|\,X\text{ is }n\text{-}C\text{-spherical}\,\}
$$
of $\mod\,R$ is contravariantly finite.
\end{cor}


{\sc Acknowledgments.}
The author would like to express his deep gratitude to Shiro Goto.
He gave the author a lot of valuable comments and careful suggestions on both the mathematics and the construction of this paper.
The author is convinced that this paper could not be in such a well-organized form without his great efforts.
The author wants to give his hearty thanks to Futoshi Hayasaka and Osamu Iyama, who gave the author helpful comments.
The author also thanks the referee for useful suggestions.


\end{document}